\documentclass[reqno,12pt]{amsart}
\usepackage{texdraw}

\usepackage[german,english]{babel}

\newtheorem{theorem}{Theorem}
\newtheorem{lemma}[theorem]{Lemma} 
\newtheorem*{remark}{Remark}

\theoremstyle{definition}

\DeclareMathOperator{\Pf}{Pf}
\DeclareMathOperator{\sgn}{sgn}

\def\({\left(}
\def\){\right)}
\def\[{\left[}
\def\]{\right]}
\def\fl#1{\left\lfloor#1\right\rfloor}
\def\cl#1{\left\lceil#1\right\rceil}

\def\qbin#1#2{\genfrac{[}{]}{0pt}{}{#1}{#2}}
\def\P{\mathcal P}
\def\dsum{\displaystyle\sum}
\begin{document}
\thispagestyle{empty}

\catcode`\@=11
\font\tenln    = line10
\font\tenlnw   = linew10

\thinlines
\newskip\Einheit \Einheit=0.6cm
\newcount\xcoord \newcount\ycoord
\newdimen\xdim \newdimen\ydim \newdimen\PfadD@cke \newdimen\Pfadd@cke
\PfadD@cke1pt \Pfadd@cke0.5pt
\def\PfadDicke#1{\PfadD@cke#1 \divide\PfadD@cke by2 \Pfadd@cke\PfadD@cke \multiply\PfadD@cke by2}
\long\def\LOOP#1\REPEAT{\def\BODY{#1}\ITERATE}
\def\ITERATE{\BODY \let\next\ITERATE \else\let\next\relax\fi \next}
\let\REPEAT=\fi
\def\Punkt{\hbox{\raise-2pt\hbox to0pt{\hss\scriptsize$\bullet$\hss}}}
\def\DuennPunkt(#1,#2){\unskip
  \raise#2 \Einheit\hbox to0pt{\hskip#1 \Einheit
          \raise-2.5pt\hbox to0pt{\hss\normalsize$\bullet$\hss}\hss}}
\def\NormalPunkt(#1,#2){\unskip
  \raise#2 \Einheit\hbox to0pt{\hskip#1 \Einheit
          \raise-3pt\hbox to0pt{\hss\large$\bullet$\hss}\hss}}
\def\DickPunkt(#1,#2){\unskip
  \raise#2 \Einheit\hbox to0pt{\hskip#1 \Einheit
          \raise-4pt\hbox to0pt{\hss\Large$\bullet$\hss}\hss}}
\def\Kreis(#1,#2){\unskip
  \raise#2 \Einheit\hbox to0pt{\hskip#1 \Einheit
          \raise-4pt\hbox to0pt{\hss\Large$\circ$\hss}\hss}}
\def\Diagonale(#1,#2)#3{\unskip\leavevmode
  \xcoord#1\relax \ycoord#2\relax
      \raise\ycoord \Einheit\hbox to0pt{\hskip\xcoord \Einheit
         \unitlength\Einheit
         \line(1,1){#3}\hss}}
\def\AntiDiagonale(#1,#2)#3{\unskip\leavevmode
  \xcoord#1\relax \ycoord#2\relax \advance\xcoord by -0.05\relax
      \raise\ycoord \Einheit\hbox to0pt{\hskip\xcoord \Einheit
         \unitlength\Einheit
         \line(1,-1){#3}\hss}}
\def\Pfad(#1,#2),#3\endPfad{\unskip\leavevmode
  \xcoord#1 \ycoord#2 \thicklines\ZeichnePfad#3\endPfad\thinlines}
\def\ZeichnePfad#1{\ifx#1\endPfad\let\next\relax
  \else\let\next\ZeichnePfad
    \ifnum#1=1
      \raise\ycoord \Einheit\hbox to0pt{\hskip\xcoord \Einheit
         \vrule height\Pfadd@cke width1 \Einheit depth\Pfadd@cke\hss}%
      \advance\xcoord by 1
    \else\ifnum#1=2
      \raise\ycoord \Einheit\hbox to0pt{\hskip\xcoord \Einheit
        \hbox{\hskip-1pt\vrule height1 \Einheit width\PfadD@cke depth0pt}\hss}%
      \advance\ycoord by 1
    \else\ifnum#1=3
      \raise\ycoord \Einheit\hbox to0pt{\hskip\xcoord \Einheit
         \unitlength\Einheit
         \line(1,1){1}\hss}
      \advance\xcoord by 1
      \advance\ycoord by 1
    \else\ifnum#1=4
      \raise\ycoord \Einheit\hbox to0pt{\hskip\xcoord \Einheit
         \unitlength\Einheit
         \line(1,-1){1}\hss}
      \advance\xcoord by 1
      \advance\ycoord by -1
    \else\ifnum#1=5
\advance\ycoord by -1\raise\ycoord \Einheit\hbox to0pt{\hskip\xcoord \Einheit
        \vrule height1 \Einheit width\PfadD@cke depth0pt\hss}%
    \fi\fi\fi\fi\fi
  \fi\next}
\def\hSSchritt{\leavevmode\raise-.4pt\hbox to0pt{\hss.\hss}\hskip.2\Einheit
  \raise-.4pt\hbox to0pt{\hss.\hss}\hskip.2\Einheit
  \raise-.4pt\hbox to0pt{\hss.\hss}\hskip.2\Einheit
  \raise-.4pt\hbox to0pt{\hss.\hss}\hskip.2\Einheit
  \raise-.4pt\hbox to0pt{\hss.\hss}\hskip.2\Einheit}
\def\vSSchritt{\vbox{\baselineskip.2\Einheit\lineskiplimit0pt
\hbox{.}\hbox{.}\hbox{.}\hbox{.}\hbox{.}}}
\def\DSSchritt{\leavevmode\raise-.4pt\hbox to0pt{%
  \hbox to0pt{\hss.\hss}\hskip.2\Einheit
  \raise.2\Einheit\hbox to0pt{\hss.\hss}\hskip.2\Einheit
  \raise.4\Einheit\hbox to0pt{\hss.\hss}\hskip.2\Einheit
  \raise.6\Einheit\hbox to0pt{\hss.\hss}\hskip.2\Einheit
  \raise.8\Einheit\hbox to0pt{\hss.\hss}\hss}}
\def\dSSchritt{\leavevmode\raise-.4pt\hbox to0pt{%
  \hbox to0pt{\hss.\hss}\hskip.2\Einheit
  \raise-.2\Einheit\hbox to0pt{\hss.\hss}\hskip.2\Einheit
  \raise-.4\Einheit\hbox to0pt{\hss.\hss}\hskip.2\Einheit
  \raise-.6\Einheit\hbox to0pt{\hss.\hss}\hskip.2\Einheit
  \raise-.8\Einheit\hbox to0pt{\hss.\hss}\hss}}
\def\SPfad(#1,#2),#3\endSPfad{\unskip\leavevmode
  \xcoord#1 \ycoord#2 \ZeichneSPfad#3\endSPfad}
\def\ZeichneSPfad#1{\ifx#1\endSPfad\let\next\relax
  \else\let\next\ZeichneSPfad
    \ifnum#1=1
      \raise\ycoord \Einheit\hbox to0pt{\hskip\xcoord \Einheit
         \hSSchritt\hss}%
      \advance\xcoord by 1
    \else\ifnum#1=2
      \raise\ycoord \Einheit\hbox to0pt{\hskip\xcoord \Einheit
        \hbox{\hskip-2pt \vSSchritt}\hss}%
      \advance\ycoord by 1
    \else\ifnum#1=3
      \raise\ycoord \Einheit\hbox to0pt{\hskip\xcoord \Einheit
         \DSSchritt\hss}
      \advance\xcoord by 1
      \advance\ycoord by 1
    \else\ifnum#1=4
      \raise\ycoord \Einheit\hbox to0pt{\hskip\xcoord \Einheit
         \dSSchritt\hss}
      \advance\xcoord by 1
      \advance\ycoord by -1
    \fi\fi\fi\fi
  \fi\next}
\def\Koordinatenachsen(#1,#2){\unskip
 \hbox to0pt{\hskip-.5pt\vrule height#2 \Einheit width.5pt depth1 \Einheit}%
 \hbox to0pt{\hskip-1 \Einheit \xcoord#1 \advance\xcoord by1
    \vrule height0.25pt width\xcoord \Einheit depth0.25pt\hss}}
\def\Koordinatenachsen(#1,#2)(#3,#4){\unskip
 \hbox to0pt{\hskip-.5pt \ycoord-#4 \advance\ycoord by1
    \vrule height#2 \Einheit width.5pt depth\ycoord \Einheit}%
 \hbox to0pt{\hskip-1 \Einheit \hskip#3\Einheit 
    \xcoord#1 \advance\xcoord by1 \advance\xcoord by-#3 
    \vrule height0.25pt width\xcoord \Einheit depth0.25pt\hss}}
\def\Gitter(#1,#2){\unskip \xcoord0 \ycoord0 \leavevmode
  \LOOP\ifnum\ycoord<#2
    \loop\ifnum\xcoord<#1
      \raise\ycoord \Einheit\hbox to0pt{\hskip\xcoord \Einheit\Punkt\hss}%
      \advance\xcoord by1
    \repeat
    \xcoord0
    \advance\ycoord by1
  \REPEAT}
\def\Gitter(#1,#2)(#3,#4){\unskip \xcoord#3 \ycoord#4 \leavevmode
  \LOOP\ifnum\ycoord<#2
    \loop\ifnum\xcoord<#1
      \raise\ycoord \Einheit\hbox to0pt{\hskip\xcoord \Einheit\Punkt\hss}%
      \advance\xcoord by1
    \repeat
    \xcoord#3
    \advance\ycoord by1
  \REPEAT}
\def\Label#1#2(#3,#4){\unskip \xdim#3 \Einheit \ydim#4 \Einheit
  \def\lo{\advance\xdim by-.5 \Einheit \advance\ydim by.5 \Einheit}%
  \def\llo{\advance\xdim by-.25cm \advance\ydim by.5 \Einheit}%
  \def\loo{\advance\xdim by-.5 \Einheit \advance\ydim by.25cm}%
  \def\o{\advance\ydim by.25cm}%
  \def\ro{\advance\xdim by.5 \Einheit \advance\ydim by.5 \Einheit}%
  \def\rro{\advance\xdim by.25cm \advance\ydim by.5 \Einheit}%
  \def\roo{\advance\xdim by.5 \Einheit \advance\ydim by.25cm}%
  \def\l{\advance\xdim by-.30cm}%
  \def\r{\advance\xdim by.30cm}%
  \def\lu{\advance\xdim by-.5 \Einheit \advance\ydim by-.6 \Einheit}%
  \def\llu{\advance\xdim by-.25cm \advance\ydim by-.6 \Einheit}%
  \def\luu{\advance\xdim by-.5 \Einheit \advance\ydim by-.30cm}%
  \def\u{\advance\ydim by-.30cm}%
  \def\ru{\advance\xdim by.5 \Einheit \advance\ydim by-.6 \Einheit}%
  \def\rru{\advance\xdim by.25cm \advance\ydim by-.6 \Einheit}%
  \def\ruu{\advance\xdim by.5 \Einheit \advance\ydim by-.30cm}%
  #1\raise\ydim\hbox to0pt{\hskip\xdim
     \vbox to0pt{\vss\hbox to0pt{\hss$#2$\hss}\vss}\hss}%
}
\catcode`\@=12

\def\ldreieck{\bsegment
  \rlvec(0.866025403784439 .5) \rlvec(0 -1)
  \rlvec(-0.866025403784439 .5)  
  \savepos(0.866025403784439 -.5)(*ex *ey)
        \esegment
  \move(*ex *ey)
        }
\def\rdreieck{\bsegment
  \rlvec(0.866025403784439 -.5) \rlvec(-0.866025403784439 -.5)  \rlvec(0 1)
  \savepos(0 -1)(*ex *ey)
        \esegment
  \move(*ex *ey)
        }
\def\rhombus{\bsegment
  \rlvec(0.866025403784439 .5) \rlvec(0.866025403784439 -.5) 
  \rlvec(-0.866025403784439 -.5)  \rlvec(0 1)        
  \rmove(0 -1)  \rlvec(-0.866025403784439 .5) 
  \savepos(0.866025403784439 -.5)(*ex *ey)
        \esegment
  \move(*ex *ey)
        }
\def\RhombusA{\bsegment
  \rlvec(0.866025403784439 .5) \rlvec(0.866025403784439 -.5) 
  \rlvec(-0.866025403784439 -.5) \rlvec(-0.866025403784439 .5) 
  \savepos(0.866025403784439 -.5)(*ex *ey)
        \esegment
  \move(*ex *ey)
        }
\def\RhombusB{\bsegment
  \rlvec(0.866025403784439 .5) \rlvec(0 -1)
  \rlvec(-0.866025403784439 -.5) \rlvec(0 1) 
  \savepos(0 -1)(*ex *ey)
        \esegment
  \move(*ex *ey)
        }
\def\RhombusC{\bsegment
  \rlvec(0.866025403784439 -.5) \rlvec(0 -1)
  \rlvec(-0.866025403784439 .5) \rlvec(0 1) 
  \savepos(0.866025403784439 -.5)(*ex *ey)
        \esegment
  \move(*ex *ey)
        }
\def\hdSchritt{\bsegment 
  \lpatt(.05 .13)
  \rlvec(0.866025403784439 -.5) 
  \savepos(0.866025403784439 -.5)(*ex *ey)
        \esegment
  \move(*ex *ey)
        }
\def\vdSchritt{\bsegment
  \lpatt(.05 .13)
  \rlvec(0 -1) 
  \savepos(0 -1)(*ex *ey)
        \esegment
  \move(*ex *ey)
        }

\def\ringerl(#1 #2){\move(#1 #2)\fcir f:0 r:.15}
\def\knoten{\bsegment \fcir f:0 r:.15 \esegment}

\def\hantel(#1 #2){\fcir f:0 r:.1 \rlvec(#1 #2) \fcir f:0 r:.1}

\newbox\Adr
\setbox\Adr\vbox{
\vskip.5cm
\centerline{Institut f\"ur Mathematik der Universit\"at Wien,}
\centerline{ Nordbergstra\ss{}e 15, A-1090 Wien, Austria.}
\centerline{E-mail: \footnotesize{\tt Theresia.Eisenkoelbl@univie.ac.at}}
}

\title[$(-1)$--enumeration 
of self--complementary plane partitions]{\boldmath $(-1)$--enumeration 
of self--complementary plane partitions}
\author{Theresia Eisenk\"olbl
\box\Adr
}
\subjclass[2000]{Primary 05A15; Secondary 05B45 52C20}
\keywords{lozenge tilings, rhombus tilings, plane partitions,
determinants, pfaffians, nonintersecting lattice paths}

\begin{abstract}
{We prove a product formula for the remaining cases of the weighted
  enumeration of 
self--complementary plane partitions contained in a given box where
adding one half of an orbit of cubes and removing the other half 
of the orbit changes the sign of the weight.
We use
nonintersecting lattice path families to express this enumeration as a
Pfaffian which can be expressed in terms of the known ordinary
enumeration of self--complementary plane partitions. 
}
\end{abstract}
\maketitle

\begin{section}{Introduction} \label{introsec}
A plane partition $P$ can be defined as a finite set of points $(i,j,k)$
with $i,j,k > 0$ and if $(i,j,k) \in P$ and $1\le i'\le i$, 
$1\le j'\le j$, $1\le k'\le k$ then $(i',j',k')\in P$. We
interpret these points as midpoints of cubes and represent a plane
partition by stacks of cubes (see Figure~\ref{sceoofi}). If we have
$i\le a$, $j\le b$ and $k\le c$ for all cubes of the plane partition,
we say that the plane partition is contained in a box with sidelengths 
$a,b,c$.

\begin{figure}
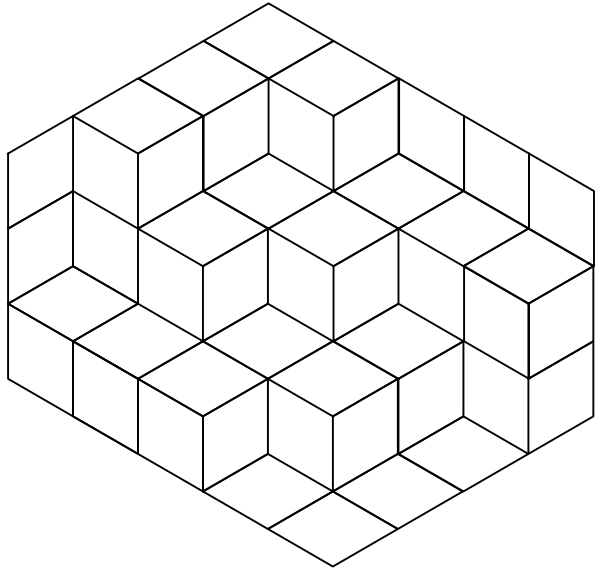

\centertexdraw{
\drawdim truecm 
\RhombusA \RhombusA \RhombusB \RhombusA \RhombusA \RhombusA \RhombusB \RhombusB
\move(-.866025 -.5)
\RhombusA \RhombusB \RhombusA \RhombusA \RhombusB \RhombusA \RhombusB\RhombusA
\move(-1.732 -1)
\RhombusA \RhombusB \RhombusA \RhombusB \RhombusA \RhombusA \RhombusB \RhombusA
\move(-2.598 -1.5)\RhombusB \RhombusB \RhombusA \RhombusA  \RhombusA \RhombusB
\RhombusA \RhombusA
\move(-1.732 -2)
\RhombusC 
\move(2.598 -.5)\RhombusC \RhombusC \RhombusC
\move(-2.598 -3.5)\RhombusC \RhombusC \RhombusC
\move(3.464 -3) \RhombusC
\move(-1.732 -4) \RhombusC 
}
\caption{A self--complementary plane partition}
\label{sceoofi}
\end{figure} 

Plane partitions were first introduced by MacMahon. One of his main results
is the following \cite[Art.~429, $x\to1$, 
proof in Art.~494]{MM}:

\smallskip
{\em The number of all plane partitions contained in a box 
with sidelengths $a,b,c$ equals}
\begin{equation} \label{box}
B(a,b,c)= \prod _{i=1} ^{a}\prod _{j=1} ^{b}\prod _{k=1} ^{c}\frac {i+j+k-1}
{i+j+k-2}= \prod _{i=1} ^{a}{\frac {(c+i)_b} {(i)_b}}, 
\end{equation}
where $(a)_n:=a(a+1)(a+2)\dots (a+n-1)$ is the rising factorial.

MacMahon also started the investigation of the number of plane
partitions with certain symmetries in a given box. These
numbers can also be expressed as product formulas similar to the one
given above. In \cite{Stan2}, Stanley introduced additional
complementation symmetries giving six new combinations of symmetries
which led to more conjectures all of which were settled in the 1980's
and 90's (see \cite{Stan2,Ku2,An3,Stem4}).

Many of these theorems come with $q$--analogs, that is, weighted
versions that record the number of cubes or orbits of cubes by a power
of $q$ and give expressions containing $q$--rising factorials instead
of rising factorials (see \cite{An1,An2,MRR2}). For plane partitions
with complementation 
symmetry, it seems to be difficult to find natural $q$--analogs. However, in
Stanley's paper a $q$--analog for self--complementary plane
partitions is given (the weight is not symmetric in the three
sidelengths, but the result is). Interestingly, upon setting $q=-1$ in
the various $q$--analogs, one consistently obtains 
enumerations of other objects, usually with additional symmetry
restraints. This observation, dubbed the ``(-1) phenomenon" has been
explained for many but not all cases 
by Stembridge (see \cite{Stem94a} and \cite{Stem94b}).

In \cite{Ku}, Kuperberg defines a $(-1)$--enumeration for all plane
partitions with complementation symmetry which admits a nice closed
product formula in almost all cases. These conjectures were solved in
Kuperberg's own paper and in the paper \cite{min} except for one case
without a nice product formula and the case of self-complementary
plane partitions in a box with some odd sidelengths which will be the
main theorem of this paper. 
We start with the precise definitions for this case.

A plane partition $P$ contained in the box $a\times b\times c$
is called
{\em self--complementary}
if  $(i,j,k) \in P \Leftrightarrow (a+1-i,b+1-j,c+1-k) \notin P$ for 
$1\le i\le a$, $1\le j\le b$, $1\le k\le c$. This means that one can
fill up the entire box by placing the plane partition and its mirror
image on top of each other. A convenient way to look at a
self--complementary plane partition is the projection to
the plane along the $(1,1,1)$--direction (see
Figure~\ref{sceoofi}). A plane partition contained in an $a\times
b\times c$--box becomes a rhombus tiling of a hexagon with
sidelengths $a,b,c,a,b,c$.
It is easy to see that self-complementary plane
partitions correspond exactly to those rhombus tilings with a $180^\circ$
rotational symmetry.

The $(-1)$--weight is defined as follows:
A self--complementary plane partition
contains exactly one half of each orbit under the operation $(i,j,k)
\mapsto (a+1-i,b+1-j,c+1-k)$.
Let a move consist of removing one half of an orbit and 
adding the other half. 
Two plane partitions are connected either by
an odd or by an even number of moves, so it is possible to define a
relative sign. The sign becomes absolute if we assign weight 1 to the
half-full plane partition (see Figure~\ref{sceoonormfi} for a box with
one side of even length and Figure~\ref{scoeenormfi} for a box with two).
\begin{figure}
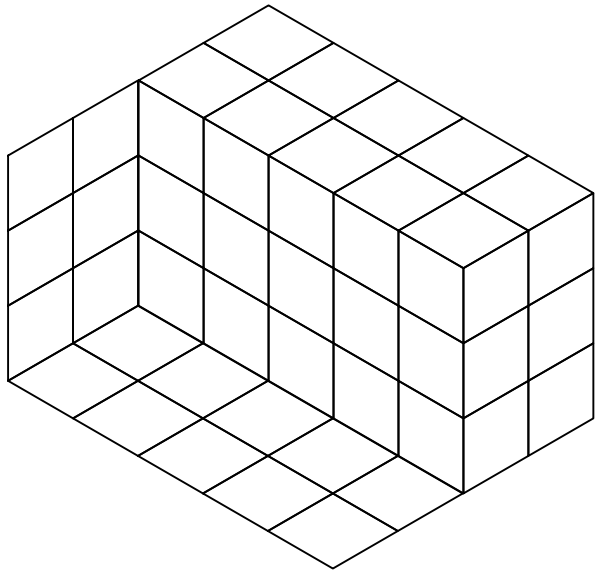

\centertexdraw{
\drawdim truecm 
\RhombusA \RhombusA \RhombusA \RhombusA \RhombusA \RhombusB \RhombusB \RhombusB
\move(-.866025 -.5)
\RhombusA \RhombusA \RhombusA \RhombusA \RhombusA \RhombusB \RhombusB\RhombusB
\move(-1.732 -1)
\RhombusB \RhombusB \RhombusB \RhombusA \RhombusA \RhombusA \RhombusA \RhombusA
\move(-2.598 -1.5)\RhombusB \RhombusB \RhombusB \RhombusA  \RhombusA \RhombusA
\RhombusA \RhombusA
\move(-.866025 -.5)
\RhombusC \RhombusC \RhombusC \RhombusC \RhombusC 
\move(-.866025 -1.5)
\RhombusC \RhombusC \RhombusC \RhombusC \RhombusC 
\move(-.866025 -2.5)
\RhombusC \RhombusC \RhombusC \RhombusC \RhombusC 
}
\caption{A plane partition of weight 1.}
\label{sceoonormfi}
\end{figure}

Therefore, this weight is
$(-1)^{n(P)}$ where $n(P)$ is the number of cubes in the ``left"
half of the box (after cutting through the sides of length $a$) if $a$
is even and $b,c$ odd or the number of cubes in the upper half of the
box (after cutting through the sides of length $b$) if $a$ is odd and
$b,c$ are even
and we
want to evaluate $\sum_{P}(-1)^{n(P)}$. For example, the plane partition in
Figure~\ref{sceoofi} has weight $(-1)^{10}=1$.

In order to be able to state the result for the $(-1)$--enumeration
more concisely, 
Stanley's result on the ordinary enumeration of self--complementary
plane partitions is needed. It will also be used as a step in the
proof of the $(-1)$--enumeration.

\begin{theorem}[Stanley \cite{Stan2}] \label{th:stanley}
The number $SC(a,b,c)$ of self--complementary plane partitions
contained in a box with sidelengths $a,b,c$ can be expressed in terms
of $B(a,b,c)$ in the following way:

\begin{align*}
B\(\tfrac a2,\tfrac b2, \tfrac c2\)^2 \quad &\text{for $a,b,c$ even,}\\
B\(\tfrac a2,\tfrac{b+1}2,\tfrac{c-1}2\) B\(\tfrac
a2,\tfrac{b-1}2,\tfrac{c+1}2\) \quad &\text{for $a$ even and $b$, $c$
  odd,}\\
B\(\tfrac {a+1}2,\tfrac{b}2,\tfrac{c}2\) B\(\tfrac
{a-1}2,\tfrac{b}2,\tfrac{c}2\)
\quad &\text{for $a$ odd and $b$, $c$
  even,}
\end{align*}
where $B(a,b,c)= \prod _{i=1} ^{a}{\frac {(c+i)_b} {(i)_b}}$ is the
number of all plane partitions in an $a\times b\times c$--box.
\end{theorem}

Note that a self-complementary plane partitions contains exactly half
of all cubes in the box. Therefore, there are no self-complementary
plane partitions in a box with three odd sidelengths.

Now we can express the $(-1)$--enumeration of self--complementary
plane partitions in terms of $SC(a,b,c)$, the ordinary
enumeration   
of self--complementary plane partitions. 

\begin{theorem} \label{main}
The enumeration of self--complementary plane
partitions in a box with sidelengths $a, b, c$ counted with weight
$(-1)^{n(P)}$ equals up to sign

\begin{align*}
B\(\frac a2,\frac b2, \frac c2\) \quad &\text{for $a,b,c$ even,}\\
SC\(\tfrac a2,\tfrac{b+1}2,\tfrac{c-1}2\) SC\(\tfrac
a2,\tfrac{b-1}2,\tfrac{c+1}2\) \quad &\text{for $a$ even and $b$, $c$
  odd}\\
SC\(\tfrac {a+1}2,\tfrac{b}2,\tfrac{c}2\) SC\(\tfrac
{a-1}2,\tfrac{b}2,\tfrac{c}2\)
\quad &\text{for $a$ odd and $b$, $c$
  even}
\end{align*}

where $SC(a,b,c)$ is given in Theorem~\ref{th:stanley} in terms of the
numbers of plane partitions contained in a box
and $n(P)$ is the number of cubes in the plane partition $P$ that are
  not in the half-full plane partition (see
Figure~\ref{sceoonormfi} and \ref{scoeenormfi}).
\end{theorem}

\begin{remark}
Note that this is zero for exactly the cases $a\equiv
2 \text{ \rm 
(mod 4)} $, $b\not\equiv c\text{ \rm (mod 4)} $ or $a$ odd, 
$b\equiv c\equiv 2 \text{ \rm (mod 4)}$ (because then the three
parameters of one factor are odd). This includes the cases where it
changes the weight if we assign 1 to another "half-full" plane partition.

Since the sides of the box play symmetric
roles this covers all cases. (For three odd sidelengths there are no
self-complementary plane partitions.) The case of three even
sidelengths has already been proved in \cite{min}.
\end{remark}

In Stanley's paper \cite{Stan2}, the theorem actually gives a
$q$--enumeration of 
plane partitions. The case $q=-1$ gives the same result as the theorem
above if one or more side has odd length, but for even sidelengths,
Stanley's theorem gives $SC\(a/2,b/2, c/2\)^2$ which
does not equal $B\( a/2, b/2,  c/2\)$. 
While the result is the same if some of the sidelengths are odd, 
the weights of individual plane partitions are different.

{\noindent \bf Outline of the proof}

{\noindent \bf Step 1: From plane partitions to families of
  nonintersecting lattice paths.}

First, we adjust a well-known bijection between 
plane partitions and families of nonintersecting lattice paths to
rephrase the problem as a path enumeration problem
(see Figure~\ref{scpathseoo} 
to get an idea).

{\noindent \bf Step 2: From lattice paths to a sum of Pfaffians}

By the main theorem on nonintersecting lattice paths (see
Lemma~\ref{gv}), this enumeration can be expressed as a sum of
determinants (see Lemma~\ref{minorsum}). 

{\noindent \bf Step 3: The sum of determinants is a single Pfaffian}

This sum can be expressed as a Pfaffian (see Lemma~\ref{pftermlemma})
by a theorem of Ishikawa and Wakayama
(see Lemma~\ref{IW}).
An analogous expression can be written down for the ordinary
enumeration of self-complementary plane partitions (see Lemma~\ref{ordenum}).

{\noindent \bf Step 4: Evaluation of the Pfaffian}

Finally, the matrix is transformed to a block matrix by elementary
row and column operations. Here, it becomes necessary to do a
case-by-case analysis according to the parity of the parameters, but
the general idea is the same in all cases.
The original entries
contain expressions with $(-1)$--binomial coefficients which are
either zero or ordinary binomial coefficients with parameters of half
the size (see \eqref{minbin}). The row and column operations involve
separating (combinations of) the even- and odd-numbered rows and
columns. Therefore, the two blocks we obtain 
have the same structure as the original matrix, but the
$(-1)$--binomial coefficients are replaced by ordinary binomial coefficients.

Now, we can identify this as certain instances of the ordinary
enumeration of self-complementary plane partitions. Since closed-form
expressions for these are already given 
by Stanley (see Theorem~\ref{th:stanley}), we can immediately derive 
the theorem.

\end{section}
\begin{section}{Proof}
{\noindent \bf Step 1: From plane partitions to families of
  nonintersecting lattice paths.}

We use the projection to the plane along the $(1,1,1)$--direction and
get immediately that self--complementary plane partitions
contained in an 
$a\times b \times c$--box are equivalent to rhombus tilings of a
hexagon with sides $a,b,c,a,b,c$ invariant under
$180^\circ$--rotation.
A tiling of this kind is clearly determined by one half of the hexagon. 

Since the sidelengths $a,b,c$ play a completely symmetric role and two
of them must have the same parity we assume without loss of generality
that $c-b$ is even and $b\le c$. The result turns 
out to be
symmetric in $b$ and $c$, so we can drop the last condition in the statement
of Theorem~\ref{main}.
Write $x$ for the positive integer $(c-b)/2$ and divide the hexagon in
half with a line parallel to the side of length $a$ (see
Figure~\ref{scpathseoo}).
As shown in the same figure, we find a bijection between these tiled
halves and families of nonintersecting lattice paths.

\newbox\orthoeoobox
\setbox\orthoeoobox\hbox{\hskip2cm
$$
\Einheit.75cm
\Gitter(7,7)(0,0)
\PfadDicke{1mm}
\Koordinatenachsen(7,7)(0,0)
\Kreis(0,3)
\Kreis(1,4)
\Kreis(2,5)
\Kreis(3,6)
\Kreis(2,1)
\Kreis(3,2)
\Kreis(5,4)
\Kreis(6,5)
\Pfad(0,3),5511\endPfad
\Pfad(1,4),1515\endPfad
\Pfad(2,5),1511\endPfad
\Pfad(3,6),1151\endPfad
\Label\lo{A_1}(0,3)
\Label\lo{A_2}(1,4)
\Label\lo{A_3}(2,5)
\Label\lo{A_4}(3,6)
\Label\ru{\kern6pt E_2}(2,1)
\Label\ru{\kern6pt E_3}(3,2)
\Label\ru{\kern6pt E_5}(5,4)
\Label\ru{\kern6pt E_6}(6,5)
\hskip-2cm
$$
}

\begin{figure}
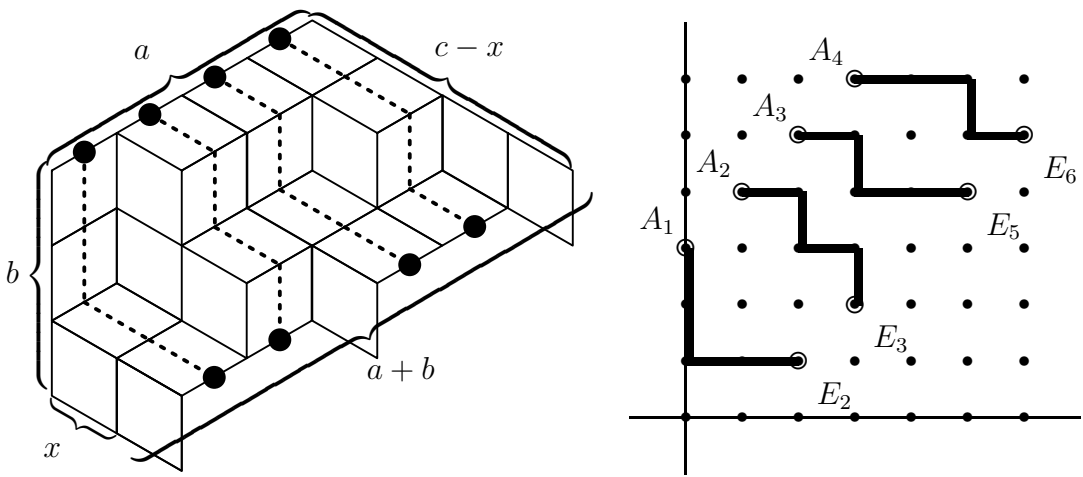

\centertexdraw{
\drawdim truecm 
\RhombusA \RhombusA \RhombusB \RhombusA 
\move(-.866025 -.5)
\RhombusA \RhombusB \RhombusA \RhombusA 
\move(-1.732 -1)
\RhombusA \RhombusB \RhombusA \RhombusB 
\move(-2.598 -1.5)\RhombusB \RhombusB \RhombusA \RhombusA  
\move(-1.732 -2)
\RhombusC 
\move(2.598 -.5)\RhombusC \RhombusC 
\move(-2.598 -3.5)\RhombusC \RhombusC 
\move(-1.732 -4) \RhombusC 
\move(.866025 -2.5) \RhombusC

\linewd.05  \lpatt(.05 .13)
\move(.433 .25) \knoten
\hdSchritt\hdSchritt\vdSchritt\hdSchritt\knoten
\move(-.433 -.25) \knoten
\hdSchritt\vdSchritt\hdSchritt\hdSchritt \knoten
\move(-1.299 -.75) \knoten \hdSchritt \vdSchritt\hdSchritt\vdSchritt \knoten
\move(-2.165 -1.25)
\knoten
\vdSchritt\vdSchritt\hdSchritt\hdSchritt \knoten

\rtext td:60 (2.6 -.6){$\left. \vbox{\vskip2.2cm}\right\}$} 
\rtext td:-60 (-2.8 -.2){$\left\{\vbox{\vskip2.2cm}\right.$} 
\htext(-3.2 -4.4){$b\left\{\vbox{\vskip1.6cm}\right.$} 
\rtext td:60 ( -2.2  -5.05) {$\left\{\vbox{\vskip.6cm}\right.$} 
\rtext td:-60 (-1.6 -3.5) {$\left.\vbox{\vskip3.7cm}\right\}$} 
\htext(2.5 0){$c-x$}
\htext(-1.5 0){$a$}
\htext(-2.7 -5.3){$x$}
\htext(1.6 -4.3){$a+b$}
\htext(11 0){\phantom{mm}}
}
\vskip-6cm
\unhbox\orthoeoobox
\caption{The paths for the self--complementary plane partition in
  Figure~\ref{sceoofi} and the orthogonal version. ($x=\frac{c-b}2$)}
\label{scpathseoo} 
\end{figure}

\newbox\orthooeebox
\setbox\orthooeebox
\hbox{\hskip2.5cm
$$
\Einheit.6cm
\Gitter(9,7)(0,0)
\PfadDicke{1mm}
\Koordinatenachsen(9,7)(0,0)
\Kreis(0,2)
\Kreis(1,3)
\Kreis(2,4)
\Kreis(3,5)
\Kreis(4,6)
\Kreis(2,0)
\Kreis(4,2)
\Kreis(5,3)
\Kreis(6,4)
\Kreis(8,6)
\Pfad(0,2),5151\endPfad
\Pfad(1,3),5111\endPfad
\Pfad(2,4),1511\endPfad
\Pfad(3,5),1115\endPfad
\Pfad(4,6),1111\endPfad
\Label\lo{A_1}(0,2)
\Label\lo{A_2}(1,3)
\Label\lo{A_3}(2,4)
\Label\lo{A_4}(3,5)
\Label\lo{A_5}(4,6)
\Label\ru{\kern6pt E_1}(2,0)
\Label\ru{\kern6pt E_3}(4,2)
\Label\ru{\kern6pt E_4}(5,3)
\Label\ru{\kern6pt E_5}(6,4)
\Label\ru{\kern6pt E_7}(8,6)
\hskip-2cm
$$}

\begin{figure}
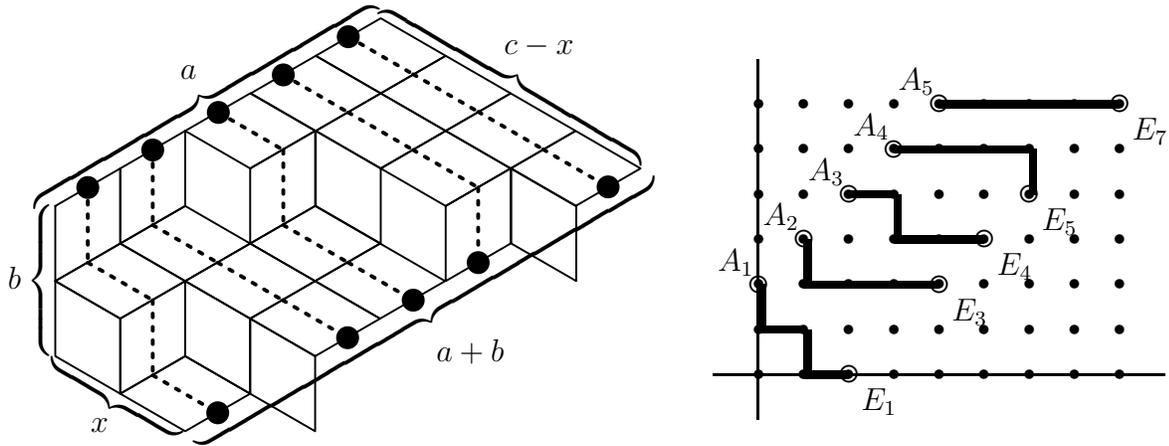

\centertexdraw{\drawdim cm
\setunitscale 1
\RhombusA\RhombusA\RhombusA\RhombusA
\move(-.866025 -.5) \RhombusA\RhombusA\RhombusA\RhombusB
\move(-1.732 -1) \RhombusA \RhombusB \RhombusA\RhombusA
\move(-2.598 -1.5)
\RhombusB\RhombusA\RhombusA\RhombusA
\move(-3.464 -2)
 \RhombusB\RhombusA
\RhombusB\RhombusA
\move(2.598 -1.5) \RhombusC
\move(0 -1) \RhombusC
\move(-1.732 -3) \RhombusC \RhombusC
\move(-3.464 -3) \RhombusC
\linewd.05  \lpatt(.05 .13)
\move(.433 .25) \knoten
\hdSchritt\hdSchritt\hdSchritt\hdSchritt\knoten
\move(-.433 -.25) \knoten
\hdSchritt\hdSchritt\hdSchritt\vdSchritt \knoten
\move(-1.299 -.75) \knoten \hdSchritt \vdSchritt\hdSchritt\hdSchritt\knoten
\move(-2.165 -1.25)
\knoten
\vdSchritt\hdSchritt\hdSchritt\hdSchritt \knoten
\move(-3.031 -1.75) \knoten
\vdSchritt\hdSchritt\vdSchritt\hdSchritt \knoten

\rtext td:60 (2.6 -.6){$\left. \vbox{\vskip2.2cm}\right\}$}
\rtext td:-60 (-3.9 -.5){$\left\{\vbox{\vskip2.8cm}\right.$}
\htext(-4.1 -4){$b\left\{\vbox{\vskip1.2cm}\right.$}
\rtext td:60 ( -2.7  -4.8) {$\left\{\vbox{\vskip1.2cm}\right.$}
\rtext td:-60 (-1.7 -3.3) {$\left.\vbox{\vskip3.7cm}\right\}$}
\htext(2.5 0){$c-x$}
\htext(-1.8 -.3){$a$}
\htext(-3 -5){$x$}
\htext(1.6 -4.1){$a+b$}
\htext(11.5 0){\phantom{mm}}
}
\vskip-5cm
\unhbox\orthooeebox
\vskip1cm
\caption{The paths for a self--complementary plane partitions with odd
  $a$. ($x=\frac{c-b}2$)}
\label{scpathsoee} 
\end{figure}

The starting points of the lattice paths are the midpoints
of the edges on the side of length $a$. The end points are the
midpoints of the edges parallel to $a$ on the opposite boundary. This
is a symmetric subset of the midpoints on the cutting line of length
$a+b$. 

The paths always follow the
rhombi of the given tiling by connecting midpoints of parallel rhombus
edges. It is easily seen that the resulting paths have no common
points (i.e. they are nonintersecting) and the tiling can be recovered
from a nonintersecting lattice path family with unit diagonal and
down steps and appropriate starting and end points. Of course, the path
families will have to be counted with the appropriate $(-1)$--weight.

\begin{figure}
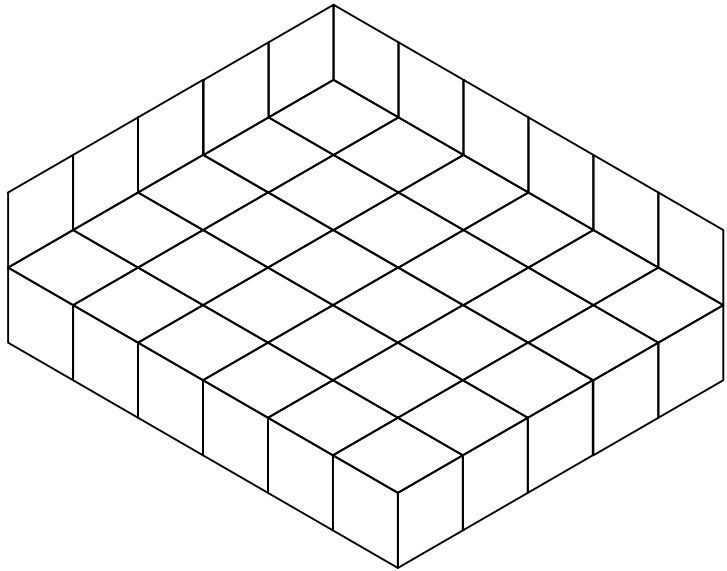

\centertexdraw{
\drawdim truecm 
\RhombusB \RhombusA \RhombusA \RhombusA \RhombusA \RhombusA \RhombusA
\RhombusB 
\move(-.866025 -.5)
\RhombusB \RhombusA \RhombusA \RhombusA \RhombusA \RhombusA \RhombusA \RhombusB 
\move(-1.732 -1)
\RhombusB \RhombusA \RhombusA \RhombusA \RhombusA \RhombusA \RhombusA \RhombusB 
\move(-2.598 -1.5)\RhombusB \RhombusA \RhombusA \RhombusA  \RhombusA \RhombusA
\RhombusA \RhombusB 
\move(-3.464 -2)\RhombusB \RhombusA \RhombusA \RhombusA  \RhombusA \RhombusA
\RhombusA \RhombusB 

\move(.866025 .5)
\RhombusC \RhombusC \RhombusC \RhombusC \RhombusC \RhombusC
\move(-3.464 -3)
\RhombusC \RhombusC \RhombusC \RhombusC \RhombusC \RhombusC
}
\caption{A plane partition of weight 1 with odd $a$.}
\label{scoeenormfi}
\end{figure}

After changing to an orthogonal coordinate system (see
Figure~\ref{scpathseoo}), 
the paths are composed of unit South and East steps
and the coordinates of the starting points are
\begin{equation} \label{start} A_i=(i-1,b+i-1)\quad \text{for
    $i=1,\dots,a$.}\end{equation} 
The end points are $a$ points chosen symmetrically among
\begin{equation} \label{end} E_j=(x+j-1,j-1)\quad \text{for
    $j=1,\dots,a+b$.}\end{equation}
Here, symmetrically means that if $E_j$ is chosen, then $E_{a+b+1-j}$
must be chosen as well.

Note that the number $a+b$ of potential end points
on the cutting line is always odd. Therefore, there is a middle one
which is either in all path families or in none according to the
parity of $a$ (see Figures~\ref{scpathseoo} and \ref{scpathsoee}).

Now the $(-1)$--weight has to be defined for the paths.
For a path from $A_i$ to $E_j$ we can use the weight
$(-1)^{\text{area}(P)}$
where area($P$) is the area between the path and the $x$--axis and
then multiply the weights of all the paths.
We have to check that the
weight changes sign if we replace a half orbit with the complementary
half orbit. If one of the affected cubes is completely inside the half 
shown in Figure~\ref{scpathseoo} or \ref{scpathsoee}, $\sum_P
\text{\rm area}(P)$ changes 
by one. If the two 
affected cubes are on the border of the figure, two symmetric
endpoints, say $E_j$ and $E_{a+b+1-j}$, are changed to $E_{j+1}$ and
$E_{a+b-j}$ or vice versa. It is easily checked
that in this case $\sum_P \text{\rm area}(P)$ changes by $j+(a+b-j)$
which is odd.

It is straightforward to check that the weight for the ``half-full"
plane partition (see Figures~\ref{sceoonormfi} and \ref{scoeenormfi})
equals
$(-1)^{a(a-2)/8}$ for $a$ even, $b,c$ odd, and 
$(-1)^{(a+b-1)c/4}$ for $a$ odd, $b,c$ even.
Therefore, we have to
multiply the path enumeration by the respective global sign.

{\noindent \bf Step 2: From lattice paths to a sum of Pfaffians}

This weight can be expressed as a product of weights on individual
steps (the exponent of $(-1)$ is just the height of the step), so the
following lemma is applicable.  
By the main theorem on nonintersecting lattice paths 
(see \cite[Lemma~1]{LindAA} or \cite[Theorem~1]{gv}) 
the weighted count of
such families of paths can be expressed as a determinant. 

\begin{lemma} \label{gv}
Let $A_1,A_2,\dots, A_n, E_1, E_2,\dots , E_n$ be integer points
meeting the following condition:
Any path from $A_i$ to $E_l$ has a common vertex
with any path from $A_j$ to $E_k$ for any $i,j,k,l$ with $i<j$ and
$k<l$.
 
Then we have
\begin{equation} \label{eq:gv}
\P ({\mathbf A} \to {\mathbf E}, \text {\rm
nonint.})=
\det_{1\le i,j \le n}{\(\P (A_i \to E_j)\)},
\end{equation}
where $\P ({A_i} \to {E_j})$ denotes the weighted enumeration of all
paths running from $A_i$ to $E_j$ and
$\P ({\mathbf A} \to {\mathbf E}, \text {\rm nonint.})$ denotes
the weighted enumeration of all families of nonintersecting lattice
paths running from $A_i$ to $E_i$ for $i=1,\dots,n$. 
\end{lemma}

The condition on the starting and end points is fulfilled in our case
because the points lie on diagonals,
so we have to find an expression for $T_{ij}=\P(A_i\to E_j)$,
the weighted enumeration of all single paths from
$A_i$ to $E_j$ in our problem. 

It is well-known that the enumeration of paths of this kind from $(x,y)$
to $(x',y')$ is given by the
$q$--binomial coefficient $\qbin {x'-x+y-y'}{x'-x}_q$ if the weight of
a path is $q^{e}$ where $e$ is the area between the
path and a horizontal line through its endpoint.

The $q$--binomial coefficient (see \cite[p. 26]{Stan1} for further
information) can be defined as
$$\qbin nk_{q}=\frac {\prod _{j=n-k+1} ^{n}(1-q^j)}
{\prod _{j=1} ^{k}(1-q^j)}.$$
Although it is not obvious from this definition, the $q$--binomial
coefficient is a polynomial in $q$. So it makes sense to put $q=-1$.

It is easy to verify that
\begin{equation}\label{minbin}
\qbin nk_{-1}
=\begin{cases} 0\quad &\text {$n$ even, $k$ odd,}\\
\binom{\fl{n/2}}{\fl{k/2}}  \quad &\text {else.} \end{cases}
\end{equation}

Taking also into account the area between horizontal line through
the endpoint and the
$x$--axis, we obtain

$$T_{ij}=\P(A_i\to E_j)=
(-1)^{(x+j-i)(j-1)}\qbin{b+x}{b+i-j}_{-1}.
$$

Now we apply Lemma~\ref{gv} to all possible sets of end points. Thus, the
$(-1)$--enumeration can be expressed as a sum of determinants which
are minors of the $a\times(a+b)$--matrix $T$:

\begin{lemma} \label{minorsum}
The $(-1)$--enumeration can be written as

\begin{multline*}
(-1)^{a(a-2)/8}
\!\!\!\!\sum _{1\le k_1<\dots<k_{a/2}\le
  (a+b-1)/2}
\!\!\!\!\!\!\!\!\!
\det(T_{k_1},\dots,T_{k_{a/2}},T_{a+b+1-k_{a/2}},\dots,T_{a+b+1-k_1})\\
\shoveright {\text{for $a$ even and $b,c$ odd,}}\\  
(-1)^{c(a+b-1)/4}\!\!\!\!\!\!\!\!
\!\!\!\!\sum _{1\le k_1<\dots<k_{(a-1)/2} \le (a+b-1)/2}
\!\!\!\!\!\!\!\!\!\!\!\!\!\!\!\!\!\!\!\!
\det(T_{k_1},\dots,T_{k_{(a-1)/2} },T_{(a+b+1)/2},
T_{a+b+1-k_{(a-1)/2}},\dots,T_{a+b+1-k_1})\\ \quad \text{for $a$ odd and
  $b,c$ even,}
\end{multline*}

where $T_{ij}$ is $(-1)^{(x+j-i)(j-1)}\qbin{b+x}{b+i-j}_{-1}$ and
$T_j$ denotes the $j$th column of $T$ which has
length $a$.
\end{lemma}

\begin{remark}
The same argument works for the ordinary enumeration, we just
have to replace $T_{ij}$ by the ordinary enumeration $\binom{b+x} {b+i-j}$.
\end{remark}

{\noindent \bf Step 3: The sum of determinants is a single Pfaffian}

Recall that the Pfaffian of a
skew--symmetric $2n\times 2n$--matrix $M$ is defined as
$$
\Pf M=\sum _{m} ^{}{\sgn m \prod _{\substack{\{i,j\}\in m\\i<j}} ^{}
{M_{ij}}},$$
where the sum runs over all
$m=\{\{m_1,m_2\},\{m_3,m_4\},\dots,\{m_{2n-1},m_{2n}\}\}$ with the conditions
$\{m_1,\dots,m_{2n}\}=\{1,\dots,2n\}$, $m_{2k-1}<m_{2k}$ and
$m_1<m_3<\dots<m_{2n-1}$. The term $\sgn m$ is the sign of the permutation
$m_1m_2m_3\dots m_{2n}$.

We will use the fact that $\(\Pf M\)^2=\det M$ and that simultaneous
row and column operations have the same effect on the Pfaffian as
ordinary row or column operations on the determinant.

Our sums of determinants can be simplified by a theorem of Ishikawa
and Wakayama \cite[Theorem 1(1)]{IshW95}  which we use 
to express the sum as a Pfaffian.
Our way of stating the theorem is taken from 
\cite[Corollary 3.2]{Okad98}. 

\begin{lemma}\label{IW}
Suppose that $n\le p$ and $n$ is even. Let $T=(t_{ik})$ be a $p\times n$
matrix and $A=(a_{kl})$ be a $p\times p$ skew-symmetric matrix.
Then we have
$$\sum _{1\le k_1<\dots<k_n\le
p}\Pf \(A_{k_1,\dots,k_n}^{k_1,\dots,k_n}\) \det(T_{k_1,\dots,k_n})=\Pf({}^tTAT),$$
where ${}^tT$ denotes the transpose of the matrix $T$, 
$T_{k_1,\dots,k_n}$ is the matrix composed of the rows of $T$ with
indices $k_1,\dots,k_n$ and $A_{k_1,\dots,k_n}^{k_1,\dots,k_n}$
is the matrix composed of the rows and columns of $A$ with
indices $k_1,\dots,k_n$.
\end{lemma}

Now specialize to $A= \begin{pmatrix}0&I_n\\ -I_n&0\end{pmatrix}$.

\begin{lemma} \label{okkor}
Let $S$ be a $2m\times 2n$--matrix with $m\le n$ and $S^{\ast}$
be the matrix $$(S_1,\dots,S_n,S_{2n},\dots,S_{n+1})$$ where $S_j$
denotes the $j$th column of $S$. Let $A$ be the matrix
$\begin{pmatrix}0&I_n\\ -I_n&0\end{pmatrix}$. Then the following
identity holds:
\begin{multline*}\sum _{1\le k_1<\dots < k_m\le n}
\det(S_{k_1},\dots,S_{k_m},S_{2n+1-k_m},\dots,S_{2n+1-k_1}) =
\Pf(S^{\ast}A({\,}^tS^{\ast}))\\
=\Pf_{1\le i,j \le 2m} \left(\sum _{k=1} ^{n} 
\left(S_{ik} S_{j,2n+1-k} - S_{jk} S_{i,2n+1-k}\right)
\right)
.
\end{multline*}
 \end{lemma}

\begin{proof}
The proof follows from Lemma~\ref{IW} with 
$A=\begin{pmatrix}0&I_n\\ -I_n&0\end{pmatrix}$ and $T={}^tS^{\ast}$. 
The sign of $\Pf \(
A_{k_1,\dots,k_m,k_1+n,\dots,k_m+n}^{k_1,\dots,k_m,k_1+n,
\dots,k_m+n}\)$ cancels exactly with the sign obtained from the 
reordering of the columns of $S$ in the determinant.
\end{proof}

Now we apply this lemma to our sums.

\begin{lemma} \label{pftermlemma}
The Pfaffians for the various $(-1)$--enumerations for $b\le c$ are
\begin{multline*}
(-1)^{a(a-2)/8} \Pf_{1\le i,j \le a} \left(\sum _{k=1} ^{\frac
      {a+b-1} {2}} \left( T_{ik} T_{j,a+b+1-k} -T_{jk} T_{i,a+b+1-k}\right)
\right),
\\
\shoveright {\text{for $a$ even and $b,c$ odd,}}\\  
\shoveleft (-1)^{c(a+b-1)/4+(a-1)/2}\Pf_{1\le i,j \le a+1} 
\left( \begin{array}{c|c}
\dsum _{k=1} ^{\frac {a+b-1} {2}} (T_{ik} T_{j,a+b+1-k} -
T_{jk}T_{i,a+b+1-k}) & T_{i,\frac {a+b+1} {2}}\\
\hline
-T_{j,\frac{a+b+1} 2} & 0 
\end{array}\right)
\\
\quad \text{for $a$ odd and
  $b,c$ even,}
\end{multline*}

where $T_{ij}=(-1)^{(x-i)(j-1)}\qbin{b+x}{b+i-j}_{-1}$ (and $x=(c-b)/2$).
\end{lemma}

\begin{proof}
In the first case, apply the lemma with $2m=a$, $2n=a+b-1$ and
$S=(T_1,\dots,T_{\frac {a+b-1} {2}}, T_{\frac {a+b+3} {2}}, \dots,
  T_{a+b})$ to obtain 
$$\Pf_{1\le i,j \le a} \left(\sum _{k=1} ^{\frac
      {a+b-1} {2}}\left( T_{ik} T_{j,a+b+1-k} -T_{jk} T_{i,a+b+1-k}\right)
\right).$$

In the second case, apply the lemma with $2m=a+1$, $2n=a+b+1$ and 
$$S=
\begin{pmatrix} 
T_1& \dots &T_{(a+b+1)/2}&
\begin{matrix}0\\0\\\vdots 
  \\0\end{matrix} & T_{(a+b+3)/2} &\dots& T_{a+b}
\\

0&\dots&0&1&0&\dots&0\end{pmatrix},$$
where the $T_j$ are columns of length $a$. 
We get 
$$\Pf_{1\le i,j \le a+1} 
\left( \begin{array}{c|c}
\sum _{k=1} ^{\frac {a+b-1} {2}} (T_{ik} T_{j,a+b+1-k} -
T_{jk}T_{i,a+b+1-k}) & T_{i,\frac {a+b+1} {2}}\\
\hline
-T_{j,\frac{a+b+1} 2} & 0 
\end{array}\right).$$

(The extra row and column correspond to an extra starting point
$A_{a+1}$ and extra end point $E_{a+b+1}$ which are only connected to
each other, so this end point must be chosen. This forces the choice of
$E_{(a+b+1)/2}$ and also gives an additional sign of $(-1)^{(a-1)/2}$). 
\end{proof}

\begin{lemma} \label{ordenum}
The Pfaffians for the ordinary enumeration $SC(a,b,c)$ for $b\le c$ are
\begin{multline*} 
\Pf_{1\le i,j\le a} \left(\sum _{k=1}^{\fl{\frac {a+b} {2}}}
\left( \binom{b+x}{b+i-k}\binom{b+x}{j+k-a-1}-
  \binom{b+x}{b+j-k}\binom{b+x}{i+k-a-1} \right) \right)
\\
\shoveright {\text{for $a$ and $c-b$ even}}\\  
\shoveleft (-1)^{(a-1)/2} \Pf_{1\le i,j \le a+1} 
\left( \begin{array}{c|c} 
\dsum _{k=1} ^{\frac {a+b-1} {2}} (\tbinom{b+x}{b+i-k}
\tbinom{b+x}{j+k-a-1}
 - \tbinom{b+x}{b+j-k} \tbinom{b+x}{k+i-a-1}) & 
\tbinom{b+x}{b+i-\frac {a+b+1} {2}}\\
\hline
-\tbinom{b+x}{b+j-\frac{a+b+1} 2 } & 0  
\end{array}\right)\\
\quad \text{for $a$ odd and
  $b,c$ even.}
\end{multline*}

\end{lemma}

\begin{proof}
Replace $T_{ij}$ by the ordinary enumeration of the respective
paths. This replaces $(-1)$--binomial coefficients by ordinary ones.
(Doing the same thing for the analogous expressions in Section~9 of \cite{min}
gives the result for the case of even sidelengths.)
\end{proof}

\begin{remark}
Of course, the closed form of these Pfaffians is known by Stanley's
theorem (see Theorem~\ref{th:stanley}). Therefore, we can use them
to evaluate the Pfaffians for the $(-1)$--enumeration.
\end{remark}

{\noindent \bf Step 4: Evaluation of the Pfaffian}

Now, the Pfaffians of Lemma~\ref{pftermlemma} can be reduced to
products of the known Pfaffians corresponding to the ordinary
enumeration. We have to do the calculations separately for different
parities of the parameters.

{\bf \boldmath Case $a,x$ even, $b,c$ odd}
We are in the first case of Lemma~\ref{pftermlemma}.
For $M_{ij}$ in $\Pf M$ we can write
\begin{multline*}
\sum _{k=1}
^{(a+b-1)/2}(-1)^{(k+1)(i+j) 
}\(\binom{(b+x-1)/2}{\fl{(b+i-k)/2}}\binom{(b+x-1)/2}{\fl{(j+k-a-1)/2}}
\right.\\
\left.-\binom{(b+x-1)/2}{\fl{(b+j-k)/2}}\binom{(b+x-1)/2}{\fl{(i+k-a-1)/2}}\)
\end{multline*}
with $1\le i,j \le a$.

Splitting the sum into terms $k=2l$ and $k=2l-1$ gives

\begin{multline} \label{matrixsplit}
\sum _{l=1}
^{\fl{(a+b-1)/4}}(-1)^{i+j}\(\binom{(b+x-1)/2}{(b-1)/2+\fl{(i-1)/2}-l+1} 
\binom{(b+x-1)/2}{\fl{(j-1)/2}+l-a/2}\right.\\
\shoveright{\left. -
\binom{(b+x-1)/2}{(b-1)/2+\fl{(j-1)/2}-l+1} 
\binom{(b+x-1)/2}{\fl{(i-1)/2}+l-a/2}
\)}\\
+\sum _{l=1}
^{\cl{(a+b-1)/4}} 
\(\binom{(b+x-1)/2}{(b-1)/2+\fl{i/2}-l+1} 
\binom{(b+x-1)/2}{\fl{j/2}+l-a/2-1}\right.\\
\left. -
\binom{(b+x-1)/2}{(b-1)/2+\fl{j/2}-l+1} 
\binom{(b+x-1)/2}{\fl{i/2}+l-a/2-1}
\)
\end{multline}

Now we apply some row and column operations to our matrix $M$. Start
with  $row(1)$, then write the differences $row(2i+1)-row(2i)$ for
$i=1,\dots,a/2-1$,  and finally $row(2i-1)+row(2i)$ 
for $i=1,\dots,a/2$. Now apply the same operations to the columns, so
that the resulting matrix is still skew--symmetric.
The new matrix has the same Pfaffian only up to sign
$(-1)^{(a/2)(a/2-1)/2}$  which cancels with the global sign in
Lemma~\ref{pftermlemma}.

Computation gives:

\begin{multline*}
M_{2i+1,j}- M_{2i,j} \\
=-\sum _{l=1} ^{\fl{(a+b-1)/4}} (-1)^j \( \binom{(b+x-1)/2+1} {(b-1)/2 +
i-l+1} \binom {(b+x-1)/2} {\fl{(j-1)/2} + l-a/2} \right. \\
\left.  - \binom {(b+x-1)/2} {(b-1)/2
+\fl{(j-1)/2} - l+1} \binom {(b+x-1)/2 +1} {i+l-a/2 } \)
\end{multline*}

Thus, apart from the first row and column, the left upper corner looks like
\begin{multline}\label{eooliob}
M_{2i+1,2j+1} - M_{2i, 2j+1} -M_{2i+1,2j} +M_{2i,2j}\\
=\sum _{l=1} ^{\fl{(a+b-1)/4}} \( \binom {(b+x-1)/2+1} {(b-1)/2 +i-l+1} 
\binom {(b+x-1)/2+1} {j+l-a/2 } \right. \\
\left. - \binom {(b+x-1)/2+1} {(b-1)/2 +j-l+1} \binom {(b+x-1)/2+1}
{i+l-a/2}  \),
\end{multline}
where $i,j=1,\dots a/2-1$.
Note how similar this is to the original matrix, only the
$(-1)$--binomial coefficients are now replaced with ordinary binomial
coefficients. The goal is to identify two blocks in the matrix which
correspond to ordinary enumeration of self--complementary plane partitions.

The right upper corner is zero (of size $(a/2-1) \times a/2$).

Furthermore,
\begin{multline*}
M_{2i-1,j}+ M_{2i,j} =
\sum _{l=1} ^{\cl{(a+b-1)/4}}  
\(   \binom {(b+x-1)/2 + 1} {(b-1)/2
+ i - l +1 } \binom {(b+x-1)/2} {\fl{j/2} +l- a/2  -1 }
\right. \\
\left.  - 
\binom{(b+x-1)/2} {(b-1)/2 +
\fl{j/2} -l +1} \binom {(b+x-1)/2 + 1} {i + l-a/2  -1 }\)
\end{multline*}

Therefore, we get for the right lower corner of the matrix

\begin{multline}\label{eooreun}
M_{2i-1,2j-1} +M_{2i,2j-1} +M_{2i-1,2j} +M_{2i,2j}\\
=\sum _{l= 1} ^{\cl{(a+b-1)/4} } \( \binom{(b+x-1)/2+1} {(b-1)/2 +i-l +1}
\binom{(b+x-1)/2+1} {j+l-a/2 -1} \right. \\
\left. -  \binom{(b+x-1)/2 +1} {(b-1)/2 +j-l +1 } 
\binom{(b+x-1)/2+1} {i+l-a/2 -1} \),
\end{multline}
where $i,j=1,\dots, a/2$.

This is almost a block matrix, only the first row and column spoil the
picture. 

Example ($a=8, b=3,c=7$):
$$\(\begin{array}{cccc|cccc}
0 &0 &1 &5 &0 &0 &-1 &-5 \\
0 &0 &3 &12 &0 &0 &0 &0 \\ 
-1 &-3 &0 &9 &0 &0 &0 &0 \\
-5 &-12 &-9 &0 &0 &0 &0 &0 \\ 
\hline
0 &0 &0 &0 &0 &-1 &-6 &-15 \\ 
0 &0 &0 &0 &1 &0 &-9 &-18 \\
1 &0 &0 &0 &6 &9 &0 &-9 \\ 
5 &0 &0 &0 &15 &18 &9 &0 
\end{array}\)$$

If {\bf \boldmath $(a/2)$ is even}, the right lower corner is an
$(a/2)\times (a/2)$--matrix  
with non-zero determinant, as we will see later, thus, we can use the
last $a/2$ rows to annihilate the second half of the first row. 
This potentially changes the entry 0 in position $(1,1)$, but
leaves everything else unchanged.
We can use the same linear  
combination on the last $a/2$ columns to annihilate the second half of
the first column. The 
resulting matrix is again skew--symmetric which means that the entry
$(1,1)$ has returned to the value 0. 
Since simultaneous row and column manipulations of this kind leave the
Pfaffian unchanged, it remains to find out
the Pfaffian of the right lower corner ($a/2\times a/2$) and the
Pfaffian of the left upper corner ($a/2\times a/2$).

The right lower block is given by Equation~\eqref{eooreun}.
This corresponds exactly to the first case of the ordinary enumeration
of self--complementary plane partitions in
Lemma~\ref{ordenum}. Therefore, the Pfaffian of this block is 
$SC(a/2,(b+1)/2,(c+1)/2)$ (which is non-zero as claimed).

The left upper $a/2\times a/2$ block (including the first row and column) is 

$$
\left(
\begin{array}{c|c}
0 &M_{1,2j+1}-M_{1,2j} \\\hline
 M_{2i+1,1}-M_{2i,1} &
\dsum _{l=1} ^{\fl{\frac{a+b-1}4}} 
\( \tbinom {(b+x+1)/2} {(b-1)/2 +i-l+1} \tbinom {(b+x+1)/2}
{j+l-a/2 }  
- \tbinom {(b+x+1)/2} {(b-1)/2 +j-l+1} \tbinom {(b+x+1)/2}
{i+l-a/2}  \)
\end{array}\right),$$

where $i,j$ run from 0 to $a/2-1$ and

\begin{align*}
M_{2i+1,1}-M_{2i,1}&= \sum _{l=1} ^{\fl{(a+b-1)/4}}  \(
\tbinom{(b+x+1)/2} {(b-1)/2 + i-l+1} \tbinom {(b+x-1)/2} {l-a/2} 
  - \tbinom {(b+x-1)/2} {(b-1)/2
- l+1} \tbinom {(b+x+1)/2 } {i+l-a/2 } \)\\
M_{1,2j+1}-M_{1,2j} &=
\sum _{l=1} ^{\fl{(a+b-1)/4}}  \( \tbinom {(b+x-1)/2} {(b-1)/2 - l+1}
\tbinom {(b+x +1)/2  } {j+l-a/2   }-
\tbinom{(b+x+1)/2} {(b-1)/2 + j-l+1} \tbinom {(b+x  -1)/2} {l-a/2} 
   \).
\end{align*}

Note that the exceptional row and column almost fit the general
pattern. We just have sometimes $(b+x-1)/2$ instead of $(b+x+1)/2$.
Replace $row(i)$ with $row(i)-row(i-1)$ for $i=1,2,\dots,a/2-1$ in
that order.
Then do the same thing for the columns.
In the resulting matrix all
occurrences of $(b+x+1)/2$ have been replaced with $(b+x-1)/2$.

After shifting the indices by one, we get
$$
 \dsum _{l=1} ^{\fl{\frac{a+b-1}4}} 
\( \tbinom {(b+x-1)/2} {(b-1)/2 +i-l} \tbinom {(b+x-1)/2} {j+l-a/2-1 }  
- \tbinom {(b+x-1)/2} {(b-1)/2 +j-l} \tbinom {(b+x-1)/2}{i+l-a/2-1}
\) ,
$$
for $i,j=1,\dots, a/2$.

The Pfaffian of this matrix can easily be identified as $SC(\frac
a2,\frac {b-1} 2,\frac {c-1} 2)$ by
Lemma~\ref{ordenum}. 
Using Theorem~\ref{th:stanley}, we obtain for the $(-1)$--enumeration
$$SC(\tfrac a 2,\tfrac {b+1} 2,\tfrac {c+1} 2)
SC(\tfrac a 2,\tfrac{b-1}2, \tfrac{c-1} 2)
=SC\(\tfrac a2,\tfrac{b+1}2,\tfrac{c-1}2\) 
SC\(\tfrac a2,\tfrac{b-1}2,\tfrac{c+1}2\),$$
which proves the main theorem in this case.

If {\bf \boldmath $(a/2)$ is odd}, we move the first row and column to
the $(a/2)$th place (which does not change the sign). 
Now we have an $(a-2)/2\times(a-2)/2$--block matrix
in the left upper corner which has non-zero determinant and thus
can be used to annihilate the first half of the exceptional row and
column similar to the previous case. By Equation~\eqref{eooliob} and
Lemma~\ref{ordenum} the Pfaffian of the left upper block is 
clearly $SC((a-2)/2,(b+1)/2,(c+1)/2)$.

For the right lower $(a+2)/2\times(a+2)/2$--block, note that the
relevant half of the exceptional column is

\begin{multline*}
M_{2i-1,1}+ M_{2i,1} =
\sum _{l=1} ^{\cl{(a+b-1)/4}}  
\(   \binom {(b+x+1)/2} {(b+1)/2
+ i - l  } \binom {(b+x-1)/2} {l- a/2  -1 }
\right. \\
\left.  - 
\binom{(b+x-1)/2} {(b+1)/2 -l } \binom {(b+x+1)/2 } {i + l-a/2  -1 }\).
\end{multline*}

We use again row and
column operations of the type $row(i)-row(i-1)$. This changes all
occurrences of $(b+x+1)/2$ to $(b+x-1)/2$ and the extra row and column
now fit the pattern in Equation~\eqref{eooreun}
with $i,j=0$. After shifting $i,j$ to $i-1,j-1$, we
identify this Pfaffian as $SC((a+2)/2,(b-1)/2,(c-1)/2)$. Again, by
Theorem~\ref{th:stanley}, the product of the two terms is exactly 
$SC(a/2,(b-1)/2,(c+1)/2) SC(a/2,(b+1)/2,(c-1)/2)$ as claimed in the theorem.

{\bf \boldmath Case $a$ even, $x$ odd, $b,c$ odd}
We start again from the first case of Lemma~\ref{pftermlemma} and have
to find the Pfaffian of the matrix

\begin{multline*}
M_{ij}=\sum _{k=1} ^{\frac {a+b-1} {2}} (-1)^{(i+j)(k+1)} 
\left( \qbin{b+x}{b+i-k}_{-1} \qbin{b+x} {j+k-a-1}_{-1} \right.\\
-
\left. \qbin{b+x} {b+j-k}_{-1} \qbin{b+x} {i+k-a-1}_{-1}\right),\quad \quad 
1\le i,j \le a.
\end{multline*}

In this case, we can simply reorder the rows and columns of the matrix
so that even indices come before odd indices. This introduces a 
sign that again cancels with $(-1)^{a(a-2)/8}$.

We have 

\begin{multline*}
M_{2i,2j-1}=\sum _{k=1} ^{\frac {a+b-1} {2}} (-1)^{k+1} 
\left( \qbin{b+x} {b+2i-k}_{-1} \qbin{b+x} {2j-2+k-a}_{-1} \right.\\
-
\left.\qbin{b+x} {b+2j-1-k}_{-1} \qbin{b+x} {2i+k-a-1}_{-1}\right).
\end{multline*}

Since $b+x$ is even and either $b+2i-k$ or $2j-2+k-a$ has to be zero,
the first product is always zero. The analogous argument for the
second product gives $M_{2i,2j-1}=M_{2j,2i-1}=0$. 

Therefore, we have 
to evaluate $(\Pf_{1\le i,j \le a/2}M_{2i,2j})(\Pf_{1\le i,j
  \le a/2}M_{2i-1,2j-1})$ which is clearly zero for {\bf \boldmath $a/2$ odd}.

Now for {\bf \boldmath $a/2$ even }
we have to evaluate the two Pfaffians.
Firstly, we substitute $k=2l-1$ to obtain for the left upper block:
\begin{multline*}
M_{2i,2j}=\sum _{l=1} ^{\cl{(a+b-1)/4}}\(
\binom{(b+x)/2}{(b+1)/2+i-l}\binom{(b+x)/2}{j-a/2+l-1}\right.\\
\left. -
\binom{(b+x)/2}{(b+1)/2+j-l}\binom{(b+x)/2}{i-a/2+l-1}\), \quad
\quad 1\le i,j \le a/2.
\end{multline*}

We can again identify the Pfaffian of this matrix as an ordinary
enumeration of self-complementary plane partitions by
Lemma~\ref{ordenum}, namely $SC(a/2,(b+1)/2,(c-1)/2)$ (here,
$(c-1)/2-(b+1)/2=(c-b)/2-1$ which is still positive because
$x=(c-b)/2$ is odd). 

Substituting $k=2l$, we obtain for the right lower block:
\begin{multline*}
M_{2i-1,2j-1}=\sum _{l=1} ^{\fl{(a+b-1)/4}}\(
\binom{(b+x)/2}{(b-1)/2+i-l}\binom{(b+x)/2}{j-1+l-a/2}\right.\\
\left. -
\binom{(b+x)/2}{(b-1)/2+j-l}\binom{(b+x)/2}{i-1+l-a/2}\), \quad
\quad 1\le i,j \le a/2.
\end{multline*}

By Lemma~\ref{ordenum} this is exactly $SC(a/2,(b-1)/2,(c+1)/2)$.\\
The product is $SC(a/2,(b+1)/2,(c-1)/2) SC(a/2,(b-1)/2,(c+1)/2)$ as
claimed in the theorem.

{\bf \boldmath Case: $a$ odd, $x$ even, $b$, $c$ even}

According to the second case of Lemma~\ref{pftermlemma} we have to
evaluate $\Pf_{1\le i,j \le a+1}M_{ij}$ for

$$
M= \left( \begin{array}{c|c}
\dsum _{k=1} ^{\frac {a+b-1} {2}} (T_{ik} T_{j,a+b+1-k} -
T_{jk}T_{i,a+b+1-k}) & T_{i,\frac {a+b+1} {2}}\\
\hline
-T_{j,\frac{a+b+1} 2} & 0 
\end{array}\right),
$$
where $T_{ij}=(-1)^{(x-i)(j-1)}\qbin{b+x}{b+i-j}_{-1}$ (and $x=(c-b)/2$).

We reorder rows and columns so that the even ones come before the odd
ones. This introduces a sign $(-1)^{(a+1)(a+3)/8}$ and gives almost a
block matrix because for $i\not= \frac {a+1} {2}$ we 
have
\begin{multline*}
M_{2i,2j-1}= \sum _{} ^{} (-1)^{k+1} 
\left(\qbin{b+x}{b+2i-k}_{-1} \qbin{b+x} {2j+k-a-2}_{-1} \right.\\-
\left.\qbin{b+x}{b+2j-k}_{-1} \qbin{b+x} {2i+k-a-2}_{-1}\right).
\end{multline*}

Since $b+x$ is even and either $b+2i-k$ or $2j+k-a-2$ is odd, we
get $M_{2i,2j-1}=0$.

Now we look at the exceptional row:

\begin{align} \notag
M_{a+1,2j}&= -T_{2j,\frac{a+b+1} 2} = 
-\qbin{b+x}{b+ 2j -\frac {a+b+1} {2} }_{-1}\\ \label{ex1}
&= \begin{cases} 
 0 \quad \quad &\text{for $\frac {a+b+1} {2}$ odd}\\ 
- \dbinom{(b+x)/2}{b/2+ j- \frac {a+b+1} {4}}
\quad \quad &\text{for $\frac {a+b+1} {2}$ even.}
\end{cases}
\end{align}

\begin{align} \notag
M_{a+1,2j-1}&= -T_{2j-1, \frac{a+b+1} 2} = 
- (-1)^{\frac {a+b+1} 2 -1}
\qbin{b+x}{b+ 2j-1 -\frac {a+b+1} {2}}_{-1}\\ \label{ex2}
&= \begin{cases} - \dbinom{(b+x)/2}{b/2+ j -1 -\frac {a+b-1} {4}}
\quad \quad &\text{for $\frac {a+b+1} {2}$ odd}\\
 0 \quad \quad &\text{for $\frac {a+b+1} {2}$ even.}
\end{cases}
\end{align}

Therefore, in the {\bf \boldmath subcase $\frac {a+b+1} {2}$ even}, 
we have a block matrix composed of two $\frac {a+1} {2} \times \frac
{a+1} {2}$--blocks. The Pfaffian is clearly zero {\bf \boldmath if
  $\frac {a+1} {2}$ 
is odd } which proves the theorem in this case.

{\bf \boldmath If $\frac {a+1} {2}$ is even,} we have two blocks.

The left upper $\frac {a+1} {2} \times \frac {a+1} {2}$--block:\\
\begin{align} \notag
M_{2i,2j}&= \sum _{k=1} ^{\frac {a+b-1} {2}} 
\left(\qbin{b+x} {b+2i-k}_{-1}  \qbin{b+x} {2j+k-a-1}_{-1} \right.\\ \notag
& \quad \quad \quad \quad \quad \quad \quad \quad \quad  \quad \quad
\quad  -\left.\qbin{b+x}{b+2j-k}_{-1}
  \qbin{b+x}{2i+k-a-1}_{-1} \right)\\ \label{2i2js}
&= \sum _{l=1} ^{\fl{\frac {a+b-1} {4}}} \left(\binom{\frac {b+x}
    {2}} {\frac {b} {2} +i-l} \binom {\frac {b+x} {2}} {j+l-\frac
    {a+1} {2}} -
\binom{\frac {b+x} {2}} {\frac {b} {2} +j-l} 
\binom {\frac {b+x} {2}} {i+l-\frac {a+1} {2}}
\right),
\end{align}
 for $i,j \not= \frac {a+1} {2}$.

We can use Equation~\eqref{ex1} and Lemma~\ref{ordenum} to see that
the left upper
Pfaffian is exactly $(-1)^{(a-3)/4} SC(\frac {a-1} {2}, \frac {b} {2}, \frac {c}
{2})$ (which is non-zero because $b/2$ is even).

The right lower block looks like 
\begin{multline}\label{2i-12j-1s} 
M_{2i-1,2j-1} = \sum _{l=1} ^{\fl{\frac {a+b+1} {4}}} \left( 
\binom{\frac {b+x}2} {\frac {b} {2} +i-l} \binom {\frac {b+x} {2}}
{j-\frac {a+1} {2} -1+l} \right.\\
-
\left.\binom{\frac {b+x}2} {\frac {b} {2} +j-l} \binom {\frac {b+x} {2}}
{i-\frac {a+1} {2} -1+l}\right),
\end{multline}

which is $SC(\frac{a+1}2, \frac {b} {2}, \frac {c} {2})$.

It can easily be checked that the signs cancel and the product of the
two terms is exactly as claimed in the theorem.

Now we look at the {\bf \boldmath subcase $\frac {a+b+1} {2}$ odd}.

Equations~\eqref{ex1} and \eqref{ex2} show that we have a block matrix
with a left upper block of size $\frac {a-1} {2}$ and a right lower
block of size $\frac {a+3} {2}$. Therefore, the Pfaffian is zero, 
{\bf \boldmath if $\frac {a-1} {2}$ is odd}, in accordance with the
claim in the theorem. 

{\bf \boldmath If $\frac {a-1} {2}$ is even},
the left upper block consists exactly
of the entries in Equation~\eqref{2i2js}. Lemma~\ref{ordenum}
identifies this Pfaffian as $SC((a-1)/2,b/2,c/2)$.
The right lower block is given by Equation~\eqref{2i-12j-1s} together
with Equation~\eqref{ex2}. We move the exceptional row and column from
the first to the last place which gives a sign change.  
By Lemma~\ref{ordenum}, the Pfaffian of this matrix
is $(-1)^{(a-1)/4} SC((a+1)/2,b/2,c/2)$.  
The signs cancel and the product of the two sub-Pfaffians is exactly
as claimed in the theorem.

{\bf \boldmath Case: $a$ odd, $x$ odd, $b$, $c$ even}

We start again with $\Pf_{1\le i,j \le a+1} M$ in the second case of
Lemma~\ref{pftermlemma}.

For $i,j < a+1$, we have 

\begin{multline*}
M_{ij}=
\dsum _{k=1} ^{\frac {a+b-1} {2}} (-1)^{(k+1)(i+j)}
\left(
\binom{(b+x-1)/2}{\fl{(b+i-k)/2}} \binom{(b+x-1)/2}{\fl{(k+j-a-1)/2}}\right.\\
- \left.\binom{(b+x-1)/2}{\fl{(b+j-k)/2}}
\binom{(b+x-1)/2}{\fl{(k+i-a-1)/2}}\right).
\end{multline*}

This is almost identical to the case of $a$ and $x$ even, therefore we
proceed similarly and split the sums for even and odd $k$. 

\begin{multline*}
M_{ij}=
\dsum _{l=1} ^{\fl{\frac {a+b-1} {4}}} (-1)^{i+j}
\left(
\binom{(b+x-1)/2}{b/2+\fl{i/2}-l}
\binom{(b+x-1)/2}{l+\fl{j/2}-(a+1)/2}\right.\\ 
\shoveright{- \left.
\binom{(b+x-1)/2}{b/2+\fl{j/2}-l}
\binom{(b+x-1)/2}{l+\fl{i/2}-(a+1)/2}\right)}\\
+ \dsum _{l=1} ^{\cl{\frac {a+b-1} {4}}} 
\left(
\binom{(b+x-1)/2}{b/2+\fl{(i+1)/2}-l} 
\binom{(b+x-1)/2}{l+\fl{(j-1)/2}-(a+1)/2}\right.\\
- \left.
\binom{(b+x-1)/2}{b/2+\fl{(j+1)/2}-l}
\binom{(b+x-1)/2}{l+\fl{(i-1)/2}-(a+1)/2}\right).
\end{multline*}

The extra row is

$$
M_{a+1,j}=-T_{j,(a+b+1)/2}=-(-1)^{(j-1)((a+b-1)/2)} 
\binom{(b+x-1)/2}{b/2+\fl{(j-(a+b+1)/2)/2}}.
$$

Now we perform the following row and column operations:
Replace the rows with $row(1)$, $row(2i+1)+row(2i)$ (for $i=1,\dots,
\frac {a-1} {2}$), $row(a+1)$, $row(2i)-row(2i-1)$ (for $i=1, \dots,
\frac{a-1}{2})$. 
Then do the same thing
for the columns. This introduces a sign of $(-1)^{(a-1)(a+5)/8}$.
All the four $(a+1)/2 \times (a+1)/2$--blocks of the new matrix have
exceptional first rows and columns. We have

\begin{multline*}
M_{2i+1,j}+ M_{2i,j}
=\dsum _{l=1} ^{\cl{\frac {a+b-1} {4}}} 
\left(
\binom{(b+x+1)/2}{b/2+i+1-l} 
\binom{(b+x-1)/2}{l+\fl{(j-1)/2}-(a+1)/2}\right.\\
- \left.
\binom{(b+x-1)/2}{b/2+\fl{(j+1)/2}-l}
\binom{(b+x+1)/2}{l+i-(a+1)/2}\right).
\end{multline*}

Therefore, the right upper block ($M_{2i+1,2j}+
M_{2i,2j}-M_{2i+1,2j-1}- M_{2i,2j-1}$) 
apart from its first row and first
column is 0.

The left upper block without its first row and column is given by
\begin{multline} \label{lastliob}
M_{2i+1,2j+1} +M_{2i,2j+1}+M_{2i+1,2j} +M_{2i,2j}\\
=\sum _{l=1} ^{\cl{\frac {a+b-1} {4}}} \left( 
\binom{(b+x+1)/2} {b/2+i+1-l} \binom{(b+x+1)/2} {l+j-(a+1)/2} \right.\\
-\left.\binom{(b+x+1)/2} {b/2+j+1-l} \binom{(b+x+1)/2} {l+i-(a+1)/2}
\right),
\end{multline}
for $i,j=1,\dots, (a-1)/2$.

The first column of the left upper block is given by

\begin{multline} \label{last1fir}
M_{2i+1,1} +M_{2i,1}
= \sum _{l=1} ^{\cl{\frac {a+b-1} {4}}} \left( 
\binom{(b+x+1)/2} {b/2+i+1-l} \binom{(b+x-1)/2} {l-(a+1)/2} \right.\\
-\left.\binom{(b+x-1)/2} {b/2+1-l} \binom{(b+x+1)/2} {l+i-(a+1)/2}\right).
\end{multline}

Furthermore, we compute

\begin{multline*}
M_{2i,j}-M_{2i-1,j}
=\sum _{l=1} ^{\fl{\frac {a+b-1} {4}}} (-1)^j
\left(\binom{(b+x+1)/2} {b/2+i-l} \binom{(b+x-1)/2}
  {l+ \fl{j/2} -(a+1)/2} \right.\\
-\left.\binom{(b+x-1)/2}{b/2+\fl{j/2}-l}
\binom{(b+x+1)/2}{l+i-(a+1)/2}\right).
\end{multline*}

The right lower block without its first row and column is given by
\begin{multline} \label{lastreun}
M_{2i,2j} -M_{2i-1,2j}-M_{2i,2j-1} +M_{2i-1,2j-1}\\
=\sum _{l=1} ^{\fl{\frac {a+b-1} {4}}} \left( 
\binom{(b+x+1)/2} {b/2+i-l} \binom{(b+x+1)/2} {l+j-(a+1)/2}\right.\\
-\left.\binom{(b+x+1)/2} {b/2+j-l} \binom{(b+x+1)/2} {l+i-(a+1)/2}
\right)
\end{multline}
for $i,j=1,\dots (a-1)/2$.

The second half of the first column is given by

\begin{multline} \label{last1sec}
M_{2i,1} -M_{2i-1,1}
= -\sum _{l=1} ^{\fl{\frac {a+b-1} {4}}} \left( 
\binom{(b+x+1)/2} {b/2+i-l} \binom{(b+x-1)/2} {l-(a+1)/2} \right.\\
-\left.\binom{(b+x-1)/2} {b/2-l} \binom{(b+x+1)/2} {l+i-(a+1)/2}\right).
\end{multline}

Finally, the other exceptional row and column are given by

\begin{multline} \label{lasta+1sec}
M_{a+1,2j}-M_{a+1,2j-1}
=-T_{2j,(a+b+1)/2}+T_{2j-1,(a+b+1)/2}\\
=-(-1)^{(a+b+1)/2-1}
\qbin{b+x}{b+2j-(a+b+1)/2}_{-1}+\qbin{b+x}{b+2j-1-(a+b+1)/2}_{-1}\\
=\begin{cases}
\binom{(b+x+1)/2}{b/2+j-(a+b+1)/4} 
 \quad \quad &\text{for $\frac {a+b+1} {2}$ even,}\\
 0 \quad \quad &\text{for $\frac {a+b+1} {2}$ odd.}
\end{cases}
\end{multline}

and 

\begin{multline}\label{lasta+1fir}
M_{2i+1,a+1}+M_{2i,a+1}
=T_{2i+1,(a+b+1)/2}+T_{2i,(a+b+1)/2}\\
=\qbin{b+x}{b+2i+1-(a+b+1)/2}_{-1}+(-1)^{(a+b+1)/2-1}\qbin{b+x}{b+2i-(a+b+1)/2}_{-1}\\
=\begin{cases}
 0 \quad \quad &\text{for $\frac {a+b+1} {2}$ even,}\\
\binom{(b+x+1)/2}{b/2+i-(a+b-1)/4} 
 \quad \quad &\text{for $\frac {a+b+1} {2}$ odd.}
\end{cases}
\end{multline}

Now we look at the {\bf \boldmath subcase $\frac {a+b+1} {2}$ even}.
Then $M_{2i+1,a+1}+M_{2i,a+1}=0$. \\
If {\bf \boldmath (a+1)/2 is even}, we will see  that the right lower
$(a+1)/2\times(a+1)/2$--matrix has a non-zero determinant and we can
treat the matrix as a block matrix despite the first exceptional row.

By simply subtracting each row and column from its successor, we can
change the left upper corner so that all $(b+x+1)/2$ become
$(b+x-1)/2$.
The first row and column now fit in with $i,j=0$ and an
index shift by one 
gives the $(a+1)/2\times (a+1)/2$--matrix:

\begin{multline*}
\sum _{l=1} ^{\cl{\frac {a+b-1} {4}}} \left( 
\binom{(b+x-1)/2} {b/2+i-l} \binom{(b+x-1)/2} {l+j-1-(a+1)/2}\right.\\
-\left.\binom{(b+x-1)/2} {b/2+j-l} \binom{(b+x-1)/2} {l+i-1-(a+1)/2}
\right).
\end{multline*}

By Lemma~\ref{ordenum}, the left upper block has Pfaffian 
$SC((a+1)/2,b/2,(c-2)/2)$. 

The right lower block is given by \eqref{lastreun} and
\eqref{lasta+1sec} with 
Pfaffian 
$(-1)^{(a-3)/4} SC((a-1)/2,b/2,(c+2)/2)$  which is not zero.

The signs cancel again, and by Theorem~\ref{th:stanley}, we have 
\begin{multline*}
SC((a+1)/2,b/2,(c-2)/2)SC((a-1)/2,b/2,(c+2)/2)\\=
SC((a-1)/2,b/2,c/2)SC((a+1)/2,b/2,c/2)
\end{multline*} as claimed in our theorem.

If {\bf \boldmath (a+1)/2 is odd}, we multiply the first row and
column by $(-1)$ and move them
to the place of the other special row and column, these are moved to
the last place. These operations change the sign.
The left upper
$(a-1)/2\times(a-1)/2$--matrix has a non-zero determinant and thus, we
can treat the matrix as a block matrix.

By Lemma~\ref{ordenum}, 
the left upper block given in \eqref{lastliob}
has Pfaffian $$SC((a-1)/2, (b+2)/2, c/2).$$ 

In the same way as in previous cases, we can use the row and column
operations of type $row(i)-row(i-1)$ to
obtain a right lower $(a+1)/2\times(a+1)/2$--block  given by

\begin{multline*}
\sum _{l=1} ^{\fl{\frac {a+b-1} {4}}} \left( 
\binom{(b+x-1)/2} {b/2+i-l-1} \binom{(b+x-1)/2} {l+j-(a+1)/2-1}\right.\\
-\left.\binom{(b+x-1)/2} {b/2+j-l-1} \binom{(b+x-1)/2}
{l+i-(a+1)/2-1}\right),
\end{multline*}

for $i,j=1,\dots,(a+1)/2$ and the row
$$
\binom{(b+x-1)/2}{b/2+j-1-(a+b+1)/4}, \quad \quad \text{ 
for $j=1,\dots,(a+1)/2)$.} $$
(Just apply the mentioned row and column operations to
$(-M_{a+1,1})= \binom{(b+x-1)/2}{b/2-(a+b+1)/4}$
and the expression
in Equation~\eqref{lasta+1sec}.)

By Lemma~\ref{ordenum} the Pfaffian of this matrix is 
$-(-1)^{(a-1)/4} SC((a+1)/2,(b-2)/2,c/2)$.

It can easily be checked that the signs cancel again and the product
of the two expressions is exactly as claimed in the theorem.

Now we look at the {\bf \boldmath subcase $\frac {a+b+1} {2}$ odd}. 

If {\bf \boldmath (a+1)/2 is odd}, the right lower $(a-1)/2\times
(a-1)/2$--block has non-zero determinant. The
second half of the $(a+1)$--row is zero by
equation~\eqref{lasta+1sec}. Therefore, we 
can eliminate the second half of the first row and column and have a
block matrix. Equation~\ref{lastreun} and Lemma~\ref{ordenum} show
that the Pfaffian of the right lower block is
$SC((a-1)/2,b/2,(c+2)/2)$.

The left upper $(a+3)/2 \times (a+3)/2$--block is given by
the expressions in \eqref{lastliob}, \eqref{last1fir},
\eqref{lasta+1fir} and $M_{a+1,1}=-\binom{(b+x-1)/2}{b/2-(a+b-1)/4}$.

The row and column operations $row(i) - row(i-1)$ 
again replace all occurrences of 
$(b+x+1)/2$ with $(b+x-1)/2$. The first row and column fits in with
$i,j=0$ (also for the $a+1$--entry), and we get an $(a+3)/2 \times
(a+3)/2$--block starting with
\begin{multline*}
\sum _{l=1} ^{\cl{(a+b-1)/4}} \left( 
\binom{(b+x-1)/2} {b/2+i-l}
\binom{(b+x-1)/2} {l+j-1-(a+1)/2}\right.\\
-
\left.\binom{(b+x-1)/2} {b/2+j-l}
\binom{(b+x-1)/2} {l+i-1-(a+1)/2}
\right),
\end{multline*}
for $i,j=1,\dots,(a+1)/2$ while the entries of the extra row are
$$-\binom{(b+x-1)/2}{b/2+j-(a+b+3)/4}.$$
Lemma~\ref{ordenum} immediately shows that the Pfaffian of this matrix
is $$(-1)^{(a-1)/4 }SC((a+1)/2, b/2, (c-2)/2).$$

The signs cancel and the product 
$$SC((a-1)/2,b/2,(c+2)/2) SC((a+1)/2,
b/2, (c-2)/2)$$ is 
equal to the expression claimed in the theorem.

If {\bf \boldmath (a+1)/2 is even}, we start by moving the first row
and column after the other special row and column. 

Now the left upper $(a+1)/2
\times (a+1)/2$--block (given by \eqref{lastliob} and
\eqref{lasta+1fir}) has non-zero determinant and can be used to
annihilate the first half of the former first row and column. By
Lemma~\ref{ordenum} the Pfaffian of the left upper block
is $$(-1)^{(a-3)/4}SC((a-1)/2,(b+2)/2,c/2).$$

The right lower block is given by \eqref{last1sec} and
\eqref{lastreun}. We multiply the first row and column by $(-1)$ and use
row and column operations similar to the previous cases to obtain the
$(a+1)/2\times (a+1)/2$--block
\begin{multline*}
\sum _{l=1} ^{\fl{\frac {a+b-1} {4}}} \left( 
\binom{(b+x-1)/2} {b/2+i-l-1} \binom{(b+x-1)/2} {l+j-1-(a+1)/2}\right.\\
-\left.\binom{(b+x-1)/2} {b/2+j-l-1} \binom{(b+x-1)/2} {l+i-1-(a+1)/2}
\right).
\end{multline*}

By Lemma~\ref{ordenum}, this is $SC((a+1)/2,(b-2)/2,c/2)$. The signs
cancel again and 
the product $SC((a-1)/2,(b+2)/2,c/2) SC((a+1)/2,(b-2)/2,c/2)$ is
easily seen to be equal to the expression in the theorem.

This case concludes the proof of the theorem. \qed

\begin{remark}
In Equation~\eqref{matrixsplit} we see that whereever $x$ occurs,
 there is actually the expression $(b+x-1)/2$.
Now replace all occurences by four different variables in the
 following way: 

\begin{multline*}
M_{ij}(m_1, m_2, n_1, n_2, a, b) \\ 
=  \sum_{l=1} ^{\fl{(a + b - 1)/4)}} 
(-1)^{i + j} \(\binom {n_1}{ (b - 1)/2 + \fl{(i - 1)/2} - l + 1}
              \binom {m_1} { -a/2 + \fl{(j - 1)/2} + l} \right.\\
           - 
\left.     \binom {n_1} {(b - 1)/2 + \fl{(j - 1)/2} - l + 1} 
              \binom {m_1} {-a/2 + \fl{(i - 1)/2} + l}\) \\
+ 
 \sum _{l= 1} ^{\cl{(a + b - 1)/4}}
\( \binom {n_2} { (b - 1)/2 + \fl{i/2} - l +1}
            \binom {m_2} {-a/2 + \fl{j/2} + l - 1}
 \right. \\
- 
 \left.         \binom {n_2} {(b - 1)/2 + \fl{j/2} - l+1} 
            \binom {m_2} {-a/2 + \fl{i/2} + l-1}\), 
\end{multline*}
 
where $i,j=1,\dots,a$.

Experimentally, the Pfaffian of this matrix is a product of linear
factors each involving only one of the four variables. Each factor
corresponds to one of the $B(r,s,t)$--factor obtained by applying
Theorem~\ref{th:stanley} to Theorem~\ref{main}.
\end{remark}

\end{section}

\end{document}